\documentclass[12pt]{amsart}
\usepackage{amsmath,amssymb,fullpage,mathdots}

\usepackage{geometry}
\usepackage{graphicx}
\usepackage{amssymb}
\usepackage{mathrsfs}
\usepackage{enumerate}
\usepackage{amsfonts}
\usepackage{amsmath}
\usepackage{cases}
\usepackage{amsfonts}
\usepackage{amsmath}
\usepackage{fancyhdr}
\title[Real zero polynomials and A. Horn's problem]{Real zero polynomials and A. Horn's problem}

\author{Lei Cao\textsuperscript{a,b}}
\address{\textsuperscript{a}School of Mathematics and Statistics, Shandong Normal University, Shandong, 250358, China;\textsuperscript{b}Department of Mathematics, Georgian Court University, Lakewood, NJ 08701, USA}
\email {lcao@georgian.edu}
\author{Hugo J. Woerdeman\textsuperscript{c}}
\address{\textsuperscript{c}Department of Mathematics, Drexel Universty, Philadelphia, PA 19104, USA}
\email {hugo@math.drexel.edu}

\thanks{HJW is partially supported by Simons Foundation grant 355645.}


\theoremstyle{plain}

\newtheorem{thm}{Theorem}[section]
\newtheorem{cor}[thm]{Corollary}

\numberwithin{equation}{section}

\newcommand{\beq}{\begin{equation}}
\newcommand{\eeq}{\end{equation}}

\theoremstyle{remark}

\numberwithin{equation}{section}

\newcommand{\bbm}{\begin{bmatrix}}
\newcommand{\ebm}{\end{bmatrix}}

\begin{document}
%

\maketitle

\begin{abstract}
A. Horn's problem concerns finding two self adjoint matrices, so that $A$, $B$, and $A+B$ have prescribed spectrum. In this paper, we show how it connects to an interpolation problem for two variable real zero polynomials and a tracial moment problem. In addition, we outline an algorithm to construct a pair $(A,B)$.
%
%
\end{abstract}

%
%

\section{Introduction}

In this paper we relate  an interpolation problem for real zero polynomials to (i) Alfred Horn's problem, and (ii) a tracial moment problem. As a constructive proof for a determinantal representation of a real zero polynomial was presented in \cite{GKVW}, it provides an outline for a constructive solution to A. Horn's problem.

A. Horn  \cite{Horn} posed the following. Given are
${\lambda}={(\lambda_1,\lambda_2,\hdots,\lambda_n)},$
${\mu}={({\mu_1},{\mu_2},\hdots,{\mu_n})},$
${\nu}={({\nu_1},{\nu_2},\hdots,{\nu_n})} \in {\mathbb R}^n$.
Under what conditions do there exist $n\times n$ Hermitian
matrices $A$ and $ B$ so that
\begin{equation} \label{mainproblem} \sigma (A) = \lambda, \sigma (B) = \mu, {\rm \ and }\ \ \sigma (A + B ) = \nu? \end{equation}
Here $\sigma (M)$ denotes the $n$ tuple of eigenvalues of the
$n\times n$ matrix $M$. For the purpose of stating the Horn inequalities we will assume
that the tuples are ordered to satisfy
\begin{equation}\label{ordered} \lambda_1 \ge \cdots \ge \lambda_n, \mu_1 \ge \cdots \ge \mu_n, \ {\rm and}\
 \nu_1 \ge \cdots \ge \nu_n. \end{equation}
Clearly, for (\ref{mainproblem}) to have a solution one needs that
\begin{equation}\label{traceeq} |\lambda|+|\mu|={\rm tr} \; A+{\rm tr} \; B={\rm tr} \;
(A+B)=|\nu|,\end{equation} where
$|(\alpha_1,\alpha_2,\ldots,\alpha_n)|:=\alpha_1+\alpha_2+\ldots+\alpha_n.$
We will refer to (\ref{traceeq}) as the {\it trace equality}. For given $\lambda, \mu$, and $\nu$, we will always assume that the trace
equality is satisfied. In \cite{Horn} A. Horn provided a
list of inequalities, all having the form,

$$\sum_{k\in K}\nu_k \leq \sum_{i\in I}\lambda_i+\sum_{j\in J}\mu_j ,$$
for certain subsets ${I, J, K}$ of $\{1,2,\ldots,n \}$ with the same cardinality $r$, with $r < n.$
Horn defined sets $T_r^n$ of triples $(I,J,K)$ by the following inductive procedure. Set

$$U^n_r = \{(I, J, K)\ | \ \sum i+\sum j= \sum k+r(r+1)/2\}.$$
When $r = 1$, set $T^n_1=U^n_1.$ Otherwise
$$T^n_r = \{(I, J, K) \in U^n_r \  | \ {\rm for\ all}\ p<r\ {\rm and \ all}\ (F,G,H)\in T^n_p, $$
$$  \sum_{f\in F} i_f+\sum_{g\in G} j_g \leq \sum_{h\in H} k_h+p(p+1)/2\}.$$
Horn conjectured that a triple $(\lambda,\mu,\nu)$ is a solution to Problem (\ref{mainproblem}) if and only if
\begin{equation}\label{Hornineq} \sum_{k\in K^c}\nu_k \geq \sum_{i\in I^c}\lambda_i+\sum_{j\in J^c}\mu_j, \end{equation}
where $I^c, J^c$ and $K^c$ are the complements of $I, J$ and $K$ in $\{1,2,\ldots,n \},$ hold for every $(I,J,K)$ in $T^n_r$ for all $r<n.$ These are called the {\it Horm inequalities}.

In 1998, Alexander Klyachko \cite{Klyachko} gave the
connection between the so-called `Saturation conjecture' and A.
Horn's conjecture. In 1999,  the Saturation Conjecture was proved
by Allen Knutson and Terence Tao \cite{KT}, and as a consequence
A. Horn's conjecture was solved. The papers \cite{Fulton} and
\cite{Stanley} provide a good overview on the history and the
solution of this problem.

Let us next make the connection to real-zero polynomials. A two variable polynomial $q(x,y)$ is called a {\it real zero polynomial} if $q$ has real coefficients and for every pair $(x,y) \in {\mathbb R}^2$ the one variable polynomial $t \mapsto q(tx,ty)$ only has real roots. Helton and Vinnikov \cite{HV} showed by algebraic geometry techniques that every real zero polynomial $q(x,y)$ with $q(0,0)=1$ may be represented as
\begin{equation}\label{rz} q(x,y) = \det (I-xA-yB )\end{equation}
where $A$ and $B$ are real symmetric matrices. More recently, in \cite{GKVW} a constructive linear algebraic proof was given of the weaker statement that representation \eqref{rz} holds with Hermitian matrices $A$ and $B$. Note that when \eqref{rz} holds, the roots of $q(x,0)$, $q(0,y)$ and $q(x,x)$ are the reciprocals of the eigenvalues of $A$, $B$ and $A+B$, respectively.  To avoid taking reciprocals and to obtain more convenient interpolation conditions we let $P=(A+B)/2$ and $R=(A-B)/2$, and consider instead
\begin{equation}\label{rzpq} p(x,y) = \det (xI -P -yR), \end{equation} which is monic in $x$. If $\sigma (A) = \lambda, \sigma (B) = \mu$, and $\sigma (A + B ) = \nu$, then
\begin{equation}\label{intpl} p(x,0) = \prod_{i=1}^n (x-\frac{\nu_i}{2} ) , p(x,1) = \prod_{i=1}^n (x-\lambda_i), p(x,-1)=\prod_{i=1}^n (x-\mu_i) . \end{equation} Thus, given $\lambda, \mu$ and $\nu$, one can attempt to find a polynomial $p(x,y)$ satisfying the interpolation conditions \eqref{intpl} and subsequently use the construction in \cite{GKVW} to find $P$ and $R$ (and thus $A$ and $B$).

The construction in \cite{GKVW} uses the fact that the roots of a polynomial are all real when a corresponding Hankel matrix of Newton sums is positive semidefinite. In our case we view  $p(x,y)=x^n + p_1(y) x^{n-1} + \cdots + p_{n-1}(y)x + p_n(y)$ as a polynomial in $x$ with coefficients depending on $y$. Note that $p_j(y)$ is a polynomial in $y$ of degree at most $j$. Now $x\mapsto p(x,y)$ has only real roots $z_1(y) , \ldots , z_n(y)$ if the Hankel matrix $H(y)=(h_{i+j}(y))_{i,j=0}^{n-1}$ of Newton sums $h_k(y) = \sum_{i=1}^n z_i(y)^k$ is positive semidefinite. The relation between the Newton sums and the coefficients of the polynomial are as follows:
\begin{equation}\label{sp}  p_1(y) =- h_1(y) , p_j(y) = -\frac1j (h_j(y) + \sum_{k=1}^{j-1} h_k(y) p_{j-k}(y)) , j\ge 2. \end{equation}
\begin{equation}\label{sph} h_0(y) = n , h_1(y) =- p_1(y) , h_j(y) = -jp_j(y) - \sum_{k=1}^{j-1} p_k(y) s_{j-k}(y) , j\ge 2. \end{equation} The transitioning between the two corresponds exactly to the Newton identities in the theory of symmetric polynomials; see, e.g., \cite{McD}.
It is important to note that $p_j (y)  \equiv 0, j \ge n+1$. In fact, in our case we just need to enforce that
\begin{equation}\label{sp0}  p_{n+1}(y) \equiv 0, p_{n+2}(y) \equiv 0 , \ldots , p_{2n-2}(y) \equiv 0 .\end{equation}
If one introduces  the companion matrix
\begin{equation}\label{comp} C(y) = \begin{pmatrix} 0 & \cdots & 0 & -p_n(y) \cr 1 & & 0 & -p_{n-1}(y) \cr & \ddots & & \vdots \cr 0 & & 1 & -p_1(y) \end{pmatrix} , \end{equation} then it is straightforward to check that \eqref{sp0} implies
\begin{equation}\label{HC} H(y) C(y) = C(y)^T H(y) . \end{equation} Since $H(y) \ge 0$ for all $y\in {\mathbb R}$, we can factor $H(y) = Q(y)^* Q(y)$ with $Q(y)$ a matrix polynomial of degree $n-1$ with $\det Q(z) \neq 0$ for ${\rm Im} \ z <0$ (see \cite{RR}, \cite{HH}, \cite[Section 2.7]{BW}). If $Q(y)$ is invertible (guaranteed when $x\mapsto p(x,y)$ has simple roots), equation \eqref{HC} implies that $M(y):= Q(y) C(y) Q(y)^{-1}$ is selfadjoint for $y\in{\mathbb R}$. Amazingly (see \cite[Section 4]{GKVW}), $M(y)$ is a degree one matrix polynomial (thus representable as $P+Ry$), yielding the desired determinantal representation \eqref{rzpq}. Indeed, by the property of companion matrices and simple determinant rules, $p(x,y)=  \det(xI-C(y)) = \det (xI-M(y)) = \det (xI-P-Ry)$. We shall formalize the connection further in Section \ref{const}. The algorithm this outlines is discussed in Section \ref{algo}.

The Newton sums of the polynomial $x\mapsto \det (xI-P-Ry)$ are ${\rm tr} (P+Ry)^k$, $k=0,1,\ldots$. Thus it is natural in our context, to consider the scenario where the coefficients of the polynomials ${\rm tr} (P+Ry)^k$, $k=0,\ldots n$, are prescribed and one wants to find the underlying matrices $P$ and $R$. In Section \ref{BMV} we will consider this so-called tracial moment problem and show its connection to A. Horn's problem. The recently solved Bessis,
Moussa, and Villania (BMV) conjecture \cite{BMV} states that if $P$ and $R$ are positive semidefinite, then the coefficients of the polynomials ${\rm tr} (P+Ry)^k$ are nonnegative; the solution to the BMV conjecture is due to Herbert Stahl \cite{Stahl}.
For more on tracial moment problems please see \cite{BK1}, \cite{BK2}, and \cite{BCKP}.

\section{Real zero polynomials and A. Horn's problem}\label{const}


For $z=(z_1, \ldots , z_n) \in{\mathbb R}^n$ we let $$S(z) = (\sum_{k=1}^n z_k^{i+j})_{i,j=0}^{n-1} = \begin{pmatrix} n & \sum_k z_k & \cdots & \sum_k z_k^{n-1} \cr
 \sum_k z_k &  \sum_k z_k^2 & \cdots &  \sum_k z_k^{n} \cr \vdots & \vdots &\iddots & \vdots \cr  \sum_k z_k^{n-1} &  \sum_k z_k^n & \cdots  & \sum_k z_k^{2n-2} \end{pmatrix}$$ be the Hankel matrix of the Newton sums of $z$. We have the following result.

\begin{thm}\label{21} Let $\lambda , \mu , \nu \in {\mathbb R}^n$ satisfying \eqref{traceeq} be given. The following are equivalent:
\begin{itemize} \item[(i)] There exist real symmetric matrices $A$ and $B$ so that $\sigma (A) = \lambda , \sigma (B) = \mu $ and $\sigma (A+B) =\nu$. \item[(ii)] There exist Hermitian matrices $A$ and $B$ so that $\sigma (A) = \lambda , \sigma (B) = \mu $ and $\sigma (A+B) =\nu$.
\item[(iii)] There exists a matrix polynomial $H(y) = \sum_{i=0}^{2n-2} H_i y^i =(h_{i+j}(y))_{i,j=0}^{n-1}$, where each $H_i=(h^{(i)}_{r+s})_{r,s=0}^{n-1}\in {\mathbb R}^{n\times n}$ is a Hankel matrix with $h^{(i)}_j=0$ for $j<i$,
$$ H(0) = S(\frac{\nu}{2}), H(1) = S(\lambda ), H(-1) = S(\mu),$$ $H(y) \ge 0$ for all $y \in {\mathbb R}$, and $p_j(y)$ constructed in \eqref{sp} satisfy \eqref{sp0}.
\item[(iv)] There exists a two variable real zero polynomial $q(x,y)$ so that $q(x,0)=\prod_{i=1}^n (1-\lambda_i x)$, $q(0,y)=\prod_{i=1}^n (1-\mu_i y)$ and $q(x,x)=\prod_{i=1}^n (1-\nu_i x)$.
\end{itemize}
\end{thm}

{\it Proof}.  (i)$\to$(ii). Trivial

(ii)$\to$(iii). Let $P=\frac12 (A+B)$ and $R=\frac12 (A-B)$, and put $p(x,y) = \det (xI_n - P - yR)$. Write \begin{equation}\label{sp1} p(x,y) = x^n + p_1(y)x^{n-1}  + \cdots + p_{n-1}(y) x + p_n(y) . \end{equation}
Define $h_j(y)$ via \eqref{sph}.
Then, as outlined in the introduction, $H(y)=(h_{r+s}(y))_{r,s=0}^{n-1}$ is the desired Hankel. See \cite[Section 4]{GKVW} for further details.

(iii)$\to$(iv). Given is the Hankel $H(y)=(h_{r+s}(y))_{r,s=0}^{n-1}$. Notice deg $ h_j(y) \le j$ for all $j$. Use equation \eqref{sp} to construct $p_j(y)$, $j=1, 2, \ldots $, and define $p(x,y)$ via \eqref{sp1}. Notice that \eqref{sp0} holds by assumption on $H(y)$ and also that deg $ p_j(y) \le j$ for all $j$. As $H(y)\ge 0$, we have that for each fixed $y\in{\mathbb R}$ the polynomial $p(x,y)$ in $x$ only has real roots. Follow the algorithm in the proof of \cite[Theorem 4.1]{GKVW}. That is, construct $C(y), Q(y), M(y)$ as in the introduction, and ultimately matrices $P$ and $R$, so that $p(x,y)  = \det (xI_n - P - yR)$. Letting now, $A = P+R $ and $B = P-R$, then $\det (I - xA - yB)$ gives the desired polynomial.

(iv)$\to$(i). Follows from \cite{HV}. \hfill $\square$

\bigskip
Note that the condition $H(y)\ge 0$ for all real $y$ in Theorem \ref{21}(iii) can be reformulated as a semidefinite programming feasibility condition (see, e.g., \cite[Section 2.7]{BW}). Let us state the resulting condition explicitely.

\begin{cor}\label{22}  Let $\lambda , \mu , \nu \in {\mathbb R}^n$ satisfying \eqref{traceeq} be given. The following are equivalent:
\begin{itemize} \item[(i)] There exist real symmetric matrices $A$ and $B$ so that $\sigma (A) = \lambda , \sigma (B) = \mu $ and $\sigma (A+B) =\nu$. \item[(ii)] There exists a positive semidefinite block matrix $G=(G_{ij})_{i,j=0}^{n-1}$ so that
\begin{equation}\label{hp} H_p:= \sum_{k=\max\{0,p-n+1\} }^{ \min\{ p,n-1\} } G_{p-k,k} , \ p=0,\ldots , n-2, \end{equation}
is a Hankel matrix $H_p =(h_{i+j}^{(p)})_{i,j=0}^{n-1}$ with $h_k^{(p)} =0$, $k<p$,
\begin{equation}\label{h0} H_0 = S(\frac{\nu}{2}) , \sum_{k=0}^{2n-2} H_k = S(\lambda ),  \sum_{k=0}^{2n-2} (-1)^k H_k = S(\mu ) , \end{equation}
and $h_j(y) =\sum_{k=0}^{2n-2} h_j^{(k)} y^k$ yields $p_j(y)$ constructed in \eqref{sp} satisfying \eqref{sp0}.
\item[(iii)] The Horn inequalities \eqref{Hornineq} hold.
\end{itemize}
\end{cor}

{\it Proof}.  (i)$\leftrightarrow$(ii). Follows from Theorem \ref{21}.

(i)$\leftrightarrow$(iii). This is due to Knutson and Tao \cite{KT}.
\hfill $\square$

It would be of interest to find a direct proof of (ii)$\leftrightarrow$(iii) in Corollary \ref{22}, as it would give a completely new proof of A. Horn's conjecture. For $n=2$ the equivalence (ii)$\leftrightarrow$(iii) can be seen as follows.
First note that \eqref{hp} and \eqref{h0} yield
\begin{equation}\label{eq1}G_{00}=H_0=\begin{pmatrix}2 & \frac{\nu_1+\nu_2}{2}\\ \frac{\nu_1+\nu_2}{2} & \frac{\nu_1^2+\nu_2^2}{4} \end{pmatrix},\end{equation}

\begin{equation}\label{eq2}G_{00}+ G_{01}+G_{10}+ G_{11}=H_0+H_1+H_2=\begin{pmatrix}2 & \lambda_1+\lambda_2\\ \lambda_1+\lambda_2 & \lambda_1^2+\lambda_2^2 \end{pmatrix},\end{equation}

\begin{equation}\label{eq3}G_{00}- G_{01}-G_{10}+G_{11}=H_0-H_1+H_2=\begin{pmatrix}2 & \mu_1+\mu_2\\ \mu_1+\mu_2 & \mu_1^2+\mu_2^2 \end{pmatrix}.\end{equation}
Consequently,
$$G_{01}+G_{10}=H_1=\begin{pmatrix}0 & \frac{\lambda_1+\lambda_2-\mu_1-\mu_2}{2} \\ \frac{\lambda_1+\lambda_2-\mu_1-\mu_2}{2} & \frac{\lambda_1^2+\lambda_2^2-\mu_1^2-\mu_2^2}{2}\end{pmatrix}$$
$$G_{11}=H_2=\begin{pmatrix}0 & 0 \\ 0 & \frac{2\lambda_1^2+2\lambda_2^2 +2\mu_1^2+2\mu_2^2-\nu_1^2-\nu_2^2 }{4}  \end{pmatrix}$$
Since $G=(G_{ij})_{i,j=0}^1$ is required to be positive semidefinite, and $G_{11}$ (and consequently $G$) has a zero diagonal entry, it makes the corresponding row and column in $G$ all zero. Taking this into account, the only possible choice for $G$ is
$$ G = \begin{pmatrix} G_{00} & G_{01} \cr G_{10} & G_{11} \end{pmatrix} = \begin{pmatrix}2 & \frac{\nu_1+\nu_2}{2}   & 0 & \frac{\lambda_1+\lambda_2-\mu_1-\mu_2}{2}  \\ \frac{\nu_1+\nu_2}{2} & \frac{\nu_1^2+\nu_2^2}{4}   & 0 &  \frac{\lambda_1^2+\lambda_2^2-\mu_1^2-\mu_2^2}{4}  \\  0 & 0  & 0 & 0 \\ \frac{\lambda_1+\lambda_2-\mu_1-\mu_2}{2}    & \frac{\lambda_1^2+\lambda_2^2-\mu_1^2-\mu_2^2}{4} &  0& \frac{2\lambda_1^2+2\lambda_2^2 +2\mu_1^2+2\mu_2^2-\nu_1^2-\nu_2^2 }{4} \end{pmatrix}   .$$
Now remove row and column 3, and take a Schur complement with respect to the $(1,1)$ entry. This gives
\begin{eqnarray}\nonumber&&\begin{pmatrix}\frac{\nu_1^2+\nu_2^2}{4}-\frac{1}{2}(\frac{(\nu_1+\nu_2)}{2})^2 & \frac{\lambda_1^2+\lambda_2^2-\mu_1^2-\mu_2^2}{4}-\frac{1}{2}(\frac{\nu_1+\nu_2}{2})(\frac{\lambda_1+\lambda_2-\mu_1-\mu_2}{2}) \\ \frac{\lambda_1^2+\lambda_2^2-\mu_1^2-\mu_2^2}{4}-\frac{1}{2}(\frac{\nu_1+\nu_2}{2})(\frac{\lambda_1+\lambda_2-\mu_1-\mu_2}{2}) &  \frac{2\lambda_1^2+2\lambda_2^2 +2\mu_1^2+2\mu_2^2-\nu_1^2-\nu_2^2 }{4}-\frac{1}{2}(\frac{\lambda_1+\lambda_2-\mu_1-\mu_2}{2})^2\end{pmatrix}\\
\nonumber &=&\begin{pmatrix}\frac{(\nu_1-\nu_2)^2}{8}& \frac{(\lambda_1-\lambda_2)^2-(\mu_1-\mu_2)^2}{8} \\  \frac{(\lambda_1-\lambda_2)^2-(\mu_1-\mu_2)^2}{8}&  \frac{(\lambda_1^2+\mu_1^2+\mu_2^2+2(\lambda_1+\mu_1)^2+2(\lambda_1+\mu_2)^2+(\mu_1-\mu_2)^2-2\nu_1^2-2\nu_2^2)}{8}\end{pmatrix},\end{eqnarray} which is required to be
a positive semidefinite matrix.

Without loss of generality, assume $\lambda_2=0.$ Denote $\nu_1$ by $x,$ then by \eqref{traceeq} we have $\nu_2=\lambda_1+\mu_1+\mu_2-x.$ Using this and multiplying the above matrix by 8, we obtain
$$A(x):={\small \begin{pmatrix}(2x-\lambda_1-\mu_1-\mu_2)^2& (\lambda_1-\lambda_2)^2-(\mu_1-\mu_2)^2\\ \\ (\lambda_1-\lambda_2)^2-(\mu_1-\mu_2)^2&  \ \ \ \ \ \begin{matrix} \lambda_1^2+\mu_1^2+\mu_2^2+2(\lambda_1+\mu_1)^2+2(\lambda_1+\mu_2)^2+ \cr (\mu_1-\mu_2)^2-2x^2-2(\lambda_1+\mu_1+\mu_2-x)^2 \end{matrix}\end{pmatrix}}. $$
Taking
its determinant we obtain
\begin{equation}\det(A(x))=-16(x-\lambda_1-\mu_1)(x-\lambda_1-\mu_2)(x-\mu_1)(x-\mu_2)\end{equation}
Using the convention \eqref{ordered} we obtain that
\begin{equation}\label{ord2} \lambda_1+\mu_1\geq \lambda_1+\mu_2 \geq \mu_2, \lambda_1+\mu_1\geq  \mu_1 \geq \mu_2, \nu_1 = x\geq \frac{\lambda_1+\mu_1+\mu_2}{2} .\end{equation}

We need to show that $A(x)$ is positive semidefinite if and only if the Horn's inequalities \eqref{Hornineq} hold, that are
\begin{equation}\nonumber  \max \{\lambda_1+\mu_2, \mu_1 \} \leq \nu_1 \leq \lambda_1+\mu_1. \end{equation} According to \eqref{ord2}, there are two possible orders among all those four roots:
\begin{enumerate}[(i)]
\item \label{item1}$\lambda_1+\mu_1\geq \lambda_1+\mu_2 \geq \mu_1 \geq \mu_2.$
Solving $A(x)\geq 0$ under this order, we have $$ \lambda_1+\mu_2 \leq x\leq \lambda_1+\mu_1\ \  {\rm and}\ \ \mu_2\leq x\leq \mu_1 .$$ However, if $\nu_1=x\leq \mu_1,$ then $\nu_2=\lambda_1+\mu_1+\mu_2-x \geq \lambda_1+\mu_2\geq x =\nu_1$ which is a contradiction.
\item \label{item2}$\lambda_1+\mu_1\geq \mu_1 \geq \lambda_1+\mu_2 \geq  \mu_2.$
Solving $A(x)\geq 0$ under this order, we have $$ \mu_1 \leq x\leq \lambda_1+\mu_1\ \  {\rm and}\ \ \mu_2\leq x\leq \lambda_1+\mu_2 .$$ However, if $\nu_1=x < \lambda_1+\mu_2,$ then $\nu_2=\lambda_1+\mu_1+\mu_2-x > \mu_1 \geq x =\nu_1$ which is a contradiction.
\end{enumerate}
In addition, the $(2,2)$ entry of $A(x)$ can be written as $$\lambda_1^2+\mu_1^2+\mu_2^2+2[(\lambda_1+\mu_1)^2-x^2]+2[(\lambda_1+\mu_2)^2-(\lambda_1+\mu_1+\mu_2-x)^2]+ (\mu_1-\mu_2)^2 $$ which is nonnegative for either case. Thus combining \eqref{item1} and \eqref{item2}, $A(x)$ is positive semidefinite if and only if
\begin{equation}\nonumber  \max \{\lambda_1+\mu_2, \mu_1 \} \leq \nu_1 \leq \lambda_1+\mu_1. \end{equation}

\section{Algorithm for finding a solution pair $(A,B)$ for A. Horn's problem}\label{algo}

Theorem \ref{21} shows that finding $A$ and $B$ satisfying \eqref{mainproblem} is equivalent to finding a real zero polynomial $q(x,y)$ with $q(x,0)$, $q(0,y)$ and $q(x,x)$ prescribed. Using the equivalence with Corollary \ref{22}(ii), we can restate it as follows.

Find polynomials $h_j(y)$, $j=1,\ldots , 2n-2$, and $p_j(y)$, $j=1, \ldots , n$  so that
\begin{itemize}\item[(a)] $H(y) = (h_{i+j}(y))_{i,j=0}^n \ge 0 , y \in {\mathbb R} ,$ where $h_0(y)\equiv n$,
\item[(b)] $\deg h_j \le j$ and $\deg p_j (y)\le j$
\item[(c)] $H(0)=S(\frac{\nu}{2})$, $H(1)=S(\lambda )$, $H(-1) = S(\mu )$, and
\item[(d)] $H(y) C(y) = C(y)^T H(y)$, $y\in {\mathbb R}$.
\end{itemize}
Finding a solution satisfying (a), (b) and (c) can be done using semidefinite programming. Indeed, this comes down to finding a positive semidefinite block matrix $G$ as in Corollary \ref{22}(ii) satisfying linear constraints (due to (b) and (c)). Constraint (d) above introduces quadratic constraints among the unknowns. Although sometimes quadratic constraints can be converted to linear constraints (for instance, by using a Schur complement trick), we have not found a way to do that with this particular constraint. However, when $n\le  3$, constraint (d) reduces to a linear constraint. Indeed, due to the three interpolation points we can determine the polynomials $h_1(y)$, $h_2(y)$, $p_1(y)$ and $p_2(y)$, as they have degree $\le 2$. The unknown polynomials do not multiply one another in this case. We have implemented the algorithm for $n=3$. The pseudo code is as follows.

 \noindent
\underline{\hspace*{165mm}}\\[1mm]
{\bf Algorithm} \ \ Solving $3\times 3$ A. Horn's problems   \\[-2mm]
\underline{\hspace*{165mm}}\\[-2mm]

\noindent
{\bf{Input.}}  Triples $\lambda , \mu , \nu \in {\mathbb R}^3$ satisfying the trace equality \eqref{traceeq}. \\[1mm]


\noindent\textbf{Step 1.} \begin{minipage}[t]{145mm} {Compute $S(\frac{\nu}{2}) , S(\lambda) , S(\mu ) $ and the polynomials $$ h_1(y)=(\sum_k \lambda_k - \frac12 \sum_k \nu_k)y + \frac12 \sum_k \nu_k,$$  $$h_2(y) = (\frac12 (\sum_k \lambda_k^2 -  \sum_k \mu_k^2)- \frac14\sum_k \nu_k^2)y^2 + \frac12 (\sum_k \lambda_k^2 +  \sum_k \mu_k^2)y + \frac14 \sum_k \nu_k^2.$$} \end{minipage} \\[-1mm]

\noindent\textbf{Step 2.} \begin{minipage}[t]{145mm} {Use semidefinite programming to determine a Hermitian $G=(G_{ij})_{i,j=0}^{2} \ge 0$ so that
\begin{itemize}\item[(i)] $ H_0:= G_{00}=S(\frac{\nu}{2}), \sum_{i,j=0}^2 G_{ij} = S(\lambda ) , \sum_{i,j=0}^2 (-1)^{i+j} G_{ij} = S(\mu) $
\item[(ii)] $ H_1:= G_{01}+G_{10} , H_2 = G_{02}+G_{11}+G_{20}, H_3=G_{12}+G_{21} , H_4=G_{22}$ are Hankel,
\item[(iii)] the $9\times 9$ matrix $G$ has the 4th, 7th and 8th columns and rows equal to 0,
\item[(iv)] $h_3(y) =\sum_{k=0}^3 h_{3k} y^k, h_4(y)=\sum_{k=0}^4 h_{4k} y^k$, where $H_3=(h_{3,i+j})_{i,j=0}^2$ and $ H_4=(h_{4,i+j})_{i,j=0}^2$,
are so that $h_4(y) - \frac43 h_3(y)h_1(y) -\frac13 h_1(y) (-h_1(y) h_2(y) - \frac12 h_1(y) (h_2(y) -h_1(y)^2)) - \frac12 h_2(y) (h_2(y)-h_1(y)^2) \equiv 0$.
\item[(v)] the quantity Im tr$(G_{12})$ is minimized.
\end{itemize}

If no such $G$ exists, then a solution $(A,B)$ does not exist, and abort. If $G$ exists, continue to Step 3.
}\end{minipage}  \\[-1mm]

\noindent\textbf{Step 3.} \begin{minipage}[t]{145mm} {Approximate $G$ with a rank 3 matrix $\hat G$ and perform a rank factorization $$ G \approx \hat G= \begin{pmatrix} Q_0^* \cr Q_1^* \cr Q_2^* \end{pmatrix} \begin{pmatrix} Q_0 & Q_1 & Q_2 \end{pmatrix}.$$ Put $Q(y) = Q_0 + Q_1 y + Q_2 y^2$. } \end{minipage} \\[-1mm]

\noindent\textbf{Step 4.} \begin{minipage}[t]{145mm} {
Let now
$$ C(1) = \begin{pmatrix} 0 & 0 & \lambda_1 \lambda_2 \lambda_3 \cr 1 & 0 & -(\lambda_1\lambda_2 + \lambda_2\lambda_3 + \lambda_3 \lambda_1) \cr 0 & 1 & \lambda_1+\lambda_2+\lambda_3 \end{pmatrix} ,$$ $$ C(-1) = \begin{pmatrix} 0 & 0 & \mu_1 \mu_2 \mu_3 \cr 1 & 0 & -(\mu_1\mu_2 + \mu_2\mu_3 + \mu_3 \mu_1)  \cr 0 & 1 & \mu_1+\mu_2+\mu_3 \end{pmatrix} ,$$

put $A=Q(1)C(1)Q(1)^{-1}$ and $B=Q(-1)C(-1)Q(-1)^{-1}$, and replace $A$ by $\frac12 (A+A^*)$ and $B$ by $\frac12 (B+B^*)$. } \end{minipage} \\[-1mm]

\vspace{2mm}
\noindent{\bf{Output.} } Hermitian matrices $A$ and $B$ satisfying \eqref{mainproblem}. \\[-2mm]
\underline{\hspace*{165mm}}\\[-3mm]
\noindent

 Some remarks are in order. First of all, for the semidefinite programming in Step 2 we used CVX, a package for specifying and solving convex programs \cite{cvx}, \cite{GB}. Item (iv) comes from the requirement that $p_4(y)\equiv 0$.  The minimization of Im tr$(G_{12})$ is to obtain the factor $Q(y)$ whose determinant is without roots in the lower half plane; see, for instance, \cite{HH} or \cite[Section 2.7]{BW}. The optimal matrix $G$ should have rank 3 in theory (we find a 4th singular value less than $10^{-11}$). The step where we replace $A$ by $\frac12 (A+A^*)$ and $B$ by $\frac12 (B+B^*)$ is theoretically not necessary. The above algorithm only works when the eigenvalues of $A$ and $B$ have multiplicity one; in other words, when $\lambda_i\neq \lambda_j$ and $\mu_i\neq \mu_j$ for $i\neq j$. Indeed, we find $A$ and $B$ as matrices similar to a companion matrix which makes them nonderogatory. The nonsingularity of $Q(1)$ and $Q(-1)$ also relies on this assumption of non-repeated eigenvalues.

We performed 100 experiments starting with random real matrices $\hat A$ and $\hat B$, and taking the eigenvalues of $\hat A+\hat{A}^*$, $\hat B + \hat{B}^*$ and $\hat A +\hat{A}^*+\hat B +\hat{B}^*$ as the input. The resulting matrices $A$ and $B$ show a maximal error in the eigenvalues of $A+B$ of size  $2.419169309320068\times 10^{-10}$. The eigenvalues of $A$ and $B$ show an error in the order of $10^{-15}$.

We hope that this algorithm is just a first step in trying to understand A. Horn's problem from a pure linear algebraic viewpoint, and that further progress will be made in the future. For the case $n>3$ we need to find a way to handle the quadratic constraints that appear. It is possible that the study of real solution to polynomial equations and its connection to semidefinite programming (see, for instance, the book \cite{AL}) may be of use.

\section{A tracial moment problem}\label{BMV}
Given $n \in {\mathbb N}$ and real numbers $s_{k,m}$, $0 \le k \le m \le 2n-2$, with $s_{0,0} = n $, we are interested in finding $n\times n$ Hermitian matrices $P$ and $R$ so that $$ \sum_{ |w|=m, {\rm na}(R)=k} {\rm tr}\ w(P,R) = s_{m,k}, 0 \le k \le m \le 2n-2 , $$ where the sum is taken over all words $w$ in two letters of length $m$ and with $R$ appearing $k$ times. For example, when $n=3$, it means we are seek $3\times 3$ Hermitian matrices $P$ and $R$ so that
$$ s_{1,0} = {\rm tr}\  P,  s_{0,1} = {\rm tr}\  R,  s_{2,0} = {\rm tr}\  P^2 ,  s_{2,1} = 2{\rm tr}\  PR ,  s_{2,2} = {\rm tr}\  R^2,  s_{3,0} = {\rm tr}\  P^3, $$ $$  s_{3,1} = 3{\rm tr}\  P^2R,  s_{3,2} = 3{\rm tr}\  PR^2 ,  s_{3,3} = {\rm tr}\  R^3,   s_{4,0} = {\rm tr}\  P^4 ,  s_{4,1} = 4{\rm tr}\  P^3R , $$ $$  s_{4,2} = 4{\rm tr}\  P^2R^2 + {\rm tr}\  PRPR ,  s_{4,3} = {\rm tr}\  PR^3 ,  s_{4,4} = {\rm tr}\  R^4 . $$
Notice that we used the general rule ${\rm tr}\  CD = {\rm tr}\  DC$, so that for instance ${\rm tr}\  PRP= {\rm tr}\  P^2R$.

\begin{thm}\label{main} Given are $n \in {\mathbb N}$ and real numbers $s_{m,k}$, $0 \le k \le m \le 2n-2$, with $s_{0,0} = n $.  There exists  $n\times n$ Hermitian matrices $P$ and $R$ so that \begin{equation}\label{skm} \sum_{ |w|=m, {\rm na}(R)=k} {\rm tr}\ w(P,R) = s_{m,k}, 0 \le k \le m \le 2n-2 , \end{equation} if and only the Hankel valued matrix polynomial
$ H(y):= (h_{r+s} (y) )_{r,s=0}^{n-1}$, where $h_m(y) = \sum_{k=0}^m s_{k,m} y^k$, $0\le m \le 2n-2$, satisfies
$$ H(y) \ge 0 , \ y \in {\mathbb R},$$ and $p_j(y)$ defined via \eqref{sp} satisfy \eqref{sp0}.
\end{thm}

{\it Proof}. First assume that $P$ and $R$ satisfying \eqref{skm} exist. Introduce the polynomial $h_m(y) = {\rm tr}\  \ (P+yR)^m$. Looking at the coefficient of $y^k$ in $h_m(y)$ one can check directly that it equals $ \sum_{ |w|=m, {\rm na}(R)=k} {\rm tr}\ \ w(P,R) = s_{m,k}$. Fixing $y\in {\mathbb R}$, one has that $P+yR$ is Hermitian, and thus
the characteristic polynomial $p(x):= \det (xI- (P+yR))$ of $P+yR$ only has real roots. The Newton sums of the roots correspond exactly to  $h_m(y) = {\rm tr}\  \ (P+yR)^m$, and since the roots are real the Hankel of the Newton sums is positive semidefinite. This gives $H(y) \ge 0$. In addition, $p_j(y)$ constructed via \eqref{sp} are exactly the coefficients of $\det (xI- (P+yR))$ and thus \eqref{sp0} is satisfied.

Conversely, given the numbers $s_{m,k}$, we build $H(y)$ and suppose it is positive semidefinite for all $y \in {\mathbb R}$. Define now $p_j(y)$ via \eqref{sp}.
Next, introduce the two variable polynomial
 \begin{equation}\label{sp1} p(x,y) = x^n + p_1(y)x^{n-1}  + \cdots + p_{n-1}(y) x + p_n(y) . \end{equation}
For a fixed $y$, we have that $H(y)$ is the Hankel of Newton sums of the polynomial, and since $H(y) \ge 0$, we get that $p(x,y)$ only has real roots (as a polynomial in $x$). Follow now the
proof of \cite[Theorem 4.1]{GKVW}. That is, construct $C(y), Q(y), M(y)$, and ultimately Hermitian matrices $P$ and $R$, so that $p(x,y)  = \det (xI_n - (P +yR))$. Then $h_m(y) = {\rm tr}\  (P+yR)^m$, $m=0,\ldots , 2n-2$, follows, and we are done.  \hfill $\square$

\section{Conclusion}
We made the connections between A. Horn's problem, a tracial moment problem and an interpolation problem for real zero polynomials. Via these connections, we showed that one may solve any one of them to get the solutions to the other two. Due to a recent constructive proof for a determinantal representation of a real zero polynomial given in \cite{GKVW}, we provided outlines for constructive solutions to A. Horn's problem in Section \ref{const} and a tracial moment problem in Section \ref{BMV} respectively. In Section \ref{algo}, we provided an implemented algorithm for finding a solution pair $(A,B)$ for A. Horn's problem in the case when $n=3.$ For the case when $n>3$ one needs to be able to handle certain quadratic constraints within a semidefinite programming context, which may be a topic of future research. It is the hope that the connections made in this paper will lead to an increased understanding of A. Horn's problem and its solution solely within a linear algebraic
 context.

\section*{Acknowledgement}
We thank the referee for careful reading and useful suggestions which helped in identifying an oversight and led to significant improvement of the paper.

\end{document}